\documentclass[10pt]{amsart}

\usepackage{amsfonts,amsmath,latexsym,amssymb,verbatim,amsbsy,amsthm}

\usepackage[top=1in, bottom=1in, left=1in, right=1in]{geometry}

\usepackage{graphicx} 
\usepackage{bmpsize} 
 
\usepackage{array}
\usepackage[labelsep=quad,indention=10pt]{caption}
\usepackage[list=true]{subcaption}
\usepackage{multirow}

\theoremstyle{plain}

\theoremstyle{definition}

\theoremstyle{remark}

\def \a {\alpha}

\def \l {\lambda}

\def \n {\nabla}

\def \th {\theta}

\def \bu {{\bf u}}

\begin{document}

\title{Addendum: 2D homogeneous solutions to the Euler equation}
\author{Xue Luo}
\address{School of Mathematics and Systems Science, Beihang University, 37 Xueyuan Road, Haidian District, Beijing, P. R. China. 100191}
\email{xluo@buaa.edu.cn}

\author{Roman Shvydkoy}
\address{Department of Mathematics, Statistics, and Computer Science
University of Illinois at Chicago
322 Science and Engineering Offices (M/C 249)
851 S. Morgan Street
Chicago, IL 60607-7045}
\email{shvydkoy@uic.edu}

\thanks{X.L.\ acknowledges the support of National Natural Science Foundation of China (11501023) and the Fundamental Research Funds for the Central Universities (YWF-16-SXXY-001). The work of R.S.\ is partially supported by NSF grants DMS-1210896, DMS-1515705.}

\date{\today}

\begin{abstract} In this addendum note we fill in the gap left in \cite{ls} in the description of 2D homogeneous solutions to the stationary Euler system with the help of the results of \cite{sd}. This gives a complete classification of all solutions. The note includes updated classification tables, and a reverse application to the results of \cite{sd}.
\end{abstract}


\maketitle

%\tableofcontents

In \cite{ls} the authors studied classification problem for homogeneous solutions of the steady Euler equations in the plane:
\begin{equation}\label{see}
\begin{split}
\bu \cdot \n \bu + \n p & = 0 \\
 \n \cdot \bu & = 0.
 \end{split}
 \end{equation}
Solutions studied are given by $\bu = \n^\perp \Psi $, where $\Psi(r,\th) = r^{\l} \psi(\th)$, in terms of polar coordinates $(r,\th)$, and $p = r^{2(\l-1)} P(\th)$. We assumed $\l >0$ to ensure the solutions have locally finite energy, however case $\l = 0$ corresponding to point vertices has also been classified. Plugging the ansatz into \eqref{see} one finds that $P$ is constant and the couple $(P,\psi)$ satisfies
\begin{equation}\label{ODE}
\begin{split}
2 (\l-1) P & =- (\l-1) (\psi')^2 +\l^2 \psi^2 + \l \psi'' \psi,\\
\psi(0) & = \psi(2\pi).
\end{split}
\end{equation}
The problem reduces to a Hamiltonian system with the help of the Bernoulli conserved quantity $B = (2P+\l^2 \psi^2+(\psi')^2)\psi^{\frac{2}{\l}-2}$, which replaces the pressure $P$ to convert \eqref{ODE} into the system
\begin{equation}\label{sys-B}
\left\{\begin{split}
		x' & = y\\
		y' & = -\lambda^2x +\frac{\lambda-1}\lambda B x^{\frac{\lambda-2}\lambda}.
\end{split}\right.
\end{equation}
where $(x,y) = (\psi,\psi')$  are the phase coordinates.
The pressure $P = - \frac{y^2}{2} - \frac{\l^2}{2} x^2 + \frac{B}{2} x^{\frac{2\l-2}{\l}}$ plays the role of Hamiltonian. In this formulation the classification problem reduces to finding all $2\pi$-periodic solutions to \eqref{sys-B}. 

The solutions are sorted into three categories:  elliptic ones corresponding to non-vanishing stream-function $x=\psi$ (streamlines are closed loops around the origin), hyperbolic-type described by $\psi$ vanishing at two or more points, and parabolic-type corresponds to vanishing at one point only or parallel shear type. In all but elliptic category the solutions have been completely classified for all values of $\l$ (we refer to \cite{ls} for classification tables in terms of life-time function). The elliptic solutions present the main technical difficulty due to subtleties related to the monotonicity of the half-period function $T$ in the elliptic region of the phase portrait of \eqref{sys-B}. Establishing monotonicity gives a precise count of distinct $2\pi$-periodic solutions for given value of $\l$ as it allows to find the number of $2\pi/n$ values within the range of the period function. This was accomplished in \cite{ls} for all values of $\l$ other than $1 < \l < \frac43$ and by duality $\frac34<\l<1$. Moreover we could show that the known sufficient conditions for monotonicity fails in the ranges above. However, the end-point values of the period, which, say in the interval $(1,4/3)$ are given by $\frac{2\pi}{\sqrt{2\l}}$ and $\pi$ clearly indicated that in order for a $2\pi$-periodic to occur the period must cross the end-points inside the elliptic region, which seemed somewhat unlikely.  

After the publication of \cite{ls} the authors learned of another work by P. Daskalopolos  and O. Savin \cite{sd} which addressed the problem of regularity and behavior at the origin of solutions to degenerate Monge-Amp\`ere equation (in the notation of \cite{sd})
\begin{equation}\label{MA}
\det D^2 u = |x|^\a, \quad \a >- 2.
\end{equation}
The Monge-Amp\`ere equation can be obtained from the Euler in 2D by taking the divergence of the momentum equation and writing it in terms of the stream-function $\Psi$.  While in general this procedure does not establish equivalence between the two, in the homogeneous case they do become equivalent, and the authors of \cite{sd} have also written down and studied the ODE \eqref{ODE}. While the aim of \cite{sd} was not to provide a complete classification of homogeneous solutions (in fact, homogeneous or asymptotically homogeneous solutions to \eqref{MA} is just one, but not the only class of solutions to \eqref{MA} described in \cite{sd}, and \cite{sd} focuses on positive values curvature, or in our terms $P>0$ which corresponds to elliptic solutions) in the range of question the authors proved in Proposition 5.1 (i) that the period function is bounded above and below by its boundary values directly without resorting to monotonicity. Thus, it excludes the possibility of $2\pi$-periodic solutions. To draw a more precise parallel, their homogeneity parameter $\a$ corresponds to our $\l = 2+\frac{\a}{2}$, considered for $\a >-2$. Proposition 5.1 (i) of \cite{sd} excludes solutions for $-2<\a<0$, which corresponds to $1 < \l <2$, which of course, includes the range $1<\l<4/3$. For the dual range $3/4< \l < 1$, the exclusion follows from the conjugacy relation between systems $0<\l<1$ and $1<\l$  established in \cite{ls},
\[
x \to x^{1/\l}(\th/\l).
\]
With these observations we now report a complete classification for elliptic solutions, and thus all homogeneous solutions. The updated  table is presented below. We refer to \cite{ls} for notation and explanation of rescaling procedures.

\begin{table}[!ht]
\centering
\begin{tabular}{|p{.6cm}<{\centering}|p{3cm}<{\centering}|p{3.5cm}<{\centering}|p{2cm}<{\centering}|p{3.5cm}<{\centering}|}
\hline	
	\multirow{2}{*}{$B$}&\multirow{2}{*}{$P$}&\multicolumn{3}{c|}{$\l$}\\
\cline{3-5}
&&$(1,2)\cup(2,\frac92]$&$2$&$(\frac92,\infty)$\\
\hline
\multirow{2}{*}{1}&$P=P_{\max}$&\multicolumn{3}{c|}{pure rotations}\\
\cline{2-5}
&$0<P<P_{\max}$&no$^*$&all$^{**}$&$\#\{(2,\sqrt{2\l})\cap\mathbb{N}\}$\\
\hline
\end{tabular}

\medskip
\begin{tabular}{|p{.6cm}<{\centering}|p{3cm}<{\centering}|p{4.7cm}<{\centering}|p{4.7cm}<{\centering}|}
\hline
&&$(0,\frac12)\cup(\frac12,1)$&$\frac12$\\
\hline
\multirow{2}{*}{-1}&$P=P_{\min}$&\multicolumn{2}{c|}{pure rotations}\\
\cline{2-4}
&$P_{\min}<P	<0$&no&all\\
\hline
\end{tabular}
\bigskip
\begin{flushleft}
\footnotesize{$^*$} ``no" means there are no solutions\\
	\footnotesize{$^{**}$} ``all" means that all solutions are $2\pi$ periodic
\end{flushleft}
		\bigskip
	\caption{\small{The number of elliptic periodic solution corresponding to all possible values of $B$, $P$ and $\lambda$ (after rescaling). Here, $P_{\max} =  \frac{1}{2\l} \left( \frac{\l-1}{\l^3} \right)^{\l-1}$, $P_{\min} = - \frac{1}{2\l} \left( \frac{1-\l}{\l^3} \right)^{\l-1}$}.}
\end{table}

In the reverse application, our results for the range $\l \in (2, \infty)$ improve upon \cite{sd} by giving a complete answer to the conditional statement of Proposition 5.1 (ii) as to the number solutions present in the range $\a >0$.  Namely, in the notation of \cite{sd}, for $0<\a \leq 5$  there are still no $2\pi$-periodic solutions, while for $\a >5$ their precise number is the cardinality of the set $(2,\sqrt{4+\a})\cap\mathbb{N}$.

%\bibliographystyle{plain}
%\bibliography{2d-homo-addendum}

\begin{thebibliography}{1}

\bibitem{sd}
Panagiota Daskalopoulos and Ovidiu Savin.
\newblock On {M}onge-{A}mp\`ere equations with homogeneous right-hand sides.
\newblock {\em Comm. Pure Appl. Math.}, 62(5):639--676, 2009.

\bibitem{ls}
Xue Luo and Roman Shvydkoy.
\newblock 2{D} homogeneous solutions to the {E}uler equation.
\newblock {\em Comm. Partial Differential Equations}, 40(9):1666--1687, 2015.

\end{thebibliography}

\end{document}